\documentclass[french, english,12pt]{amsart}

\usepackage{babel}
\usepackage[T1]{fontenc}
\usepackage{amsmath,amssymb}

\newcommand{\lgw}{\longrightarrow}
\newcommand{\lgm}{\longmapsto}
\newcommand{\ovl}{\overline}
\newcommand{\ord}{\text{ord}}
\newcommand{\LL}{\mathcal{L}}
\newcommand{\wdh}{\widehat}

\newcommand{\la}{\lambda}
\newcommand{\OO}{\mathcal{O}}

\newcommand{\m}{\mathfrak{m}}
\newcommand{\Z}{\mathbb{Z}}

\renewcommand{\k}{\Bbbk}
\newcommand{\K}{\mathbb{K}}
\newcommand{\N}{\mathbb{N}}

\newcommand{\R}{\mathbb{R}}

\renewcommand{\b}{\beta}
\newcommand{\g}{\gamma}

\newtheorem{definition}{D\'efinition }[section]

\newtheorem{lemme}[definition]
{Lemme }
\newtheorem{theorem}[definition]
{Th\'eor\`eme }
\newtheorem{example}[definition]
{Exemple }

\newtheorem{proposition}[definition]{Proposition}

\newtheorem{remark}{Remarque}

\begin{document}

\title[Sur la lin\'earit\'e de la fonction de Artin]{Sur la
lin\'earit\'e de la fonction de Artin\\
\vspace{0.8cm}
About the linearity of the Artin function}
\author{Guillaume Rond}

\address{
Department of Mathematics \\
University of Toronto\\
Toronto, Ontario M5S 2E4\\
 Canada\\
(\email{rond@picard.ups-tlse.fr})}

\maketitle
\selectlanguage{french}
\begin{abstract}
Nous donnons ici un contre-exemple \`a une vieille conjecture en
th\'eorie des singularit\'es. Cette conjecture est que la fonction qui
apparait dans la version forte du th\'eor\`eme d'approximation de Artin
est major\'ee par une fonction affine. Tout d'abord nous faisons une
\'etude de l'approximation diophantienne entre le corps des s\'eries en
plusieurs variables et son compl\'et\'e pour la topologie $\m$-adique.
Nous montrons, \`a l'aide d'un exemple, qu'il n'existe pas de version du
th\'eor\`eme de Liouville dans ce contexte. Ce m\^eme exemple nous
fournit notre contre-exemple (th\'eor\`eme \ref{cont}). Nous appliquons
cela pour donner une nouvelle preuve du fait qu'il n'existe pas de
th\'eorie de l'\'elimination des quantificateurs pour le corps des
s\'eries en plusieurs variables.
\end{abstract}
\selectlanguage{english}
\begin{abstract}
We give here a counter-example to an old conjecture in the theory of
singularities. This conjecture is that the function that appears in the
strong Artin approximation theorem is bounded by a affine function.
First we study diophantine approximation between the field  of power
series in several variables and its completion for the $\m$-adic
topology. We show, with an example, that there is no Liouville theorem
in this case. This example give us our counter-example (cf. th\'eor\`eme
\ref{cont}). As an application, we give a new proof of the fact that
there is no theory of elimination of quantifiers for the field of
fractions of the ring of power series in several variables.
\end{abstract}

\selectlanguage{french}

\section{Introduction}

En 1969, M. Artin a prouv\'e le th\'eor\`eme suivant (version forte du
th\'eor\`eme d'approximation  de Artin) :\\
\selectlanguage{french}
\begin{theorem}\label{Ar}\cite{Ar}
Soit $I$ un id\'eal de $\OO_N[X_1,...,\,X_n]$ engendr\'e par les $f_l$
pour $1\leq l\leq r$. Il existe une fonction $\b\,:\,\N\lgw\N$ telle que :
$$\forall i\in\N\ \forall x\in \OO_N^n\text{ tels que } $$
$$\begin{array}{ccc}\left(\forall l \ \ f_l(x)\in \m_N^{\b(i)+1}\OO_N^r\
\right) &\Longrightarrow &
\left(\exists \,\ovl{x}\in \OO_N^n \text{ tel que }
\begin{array}{c}\ovl{x}-x\in\m_N^{i+1}\OO_N^n\\
\text{ et }\forall l \ \ f_l(\ovl{x})=0\ \end{array}\right)\end{array}$$
o\`u $\m_N$ est l'id\'eal maximal de $\OO_N$, l'anneau des s\'eries
formelles en $N$ variables sur un corps $\k$ quelconque.
\end{theorem}
Pour $I$  un id\'eal de $\OO_N[X_1,...,\,X_n]$, nous appellerons fonction
de Artin de $I$ la plus petite fonction $\b\,:\,\N\lgw\N$ qui v\'erifie
le th\'eor\`eme ci-dessus. Ce th\'eor\`eme nous dit donc que les
solutions approch\'ees des $f_l$ sont proches de vraies solutions pour
la topologie $\m$-adique, la fonction $\b$ \'etant une "mesure" de cette
approximation.
Dans le cas $N=1$, M. J. Greenberg \cite{Gr} a montr\'e que la fonction
de Artin d'un id\'eal  \'etait major\'ee par une fonction affine. Le
fait que dans ce cas la fonction de Artin d'un id\'eal soit born\'ee par
une fonction affine peut s'interpr\'eter en terme d'in\'egalit\'e de
type {\L}ojasiewicz sur l'espace des arcs (cf. remarque \ref{Loj}).\\
D'autre part, si $\OO_N\lgw\OO_N[X_1,...,\,X_n]/I$ est lisse et $N$
quelconque, alors la fonction de Artin de $I$ est \'egale \`a
l'identit\'e. La fonction de Artin d'un tel id\'eal peut donc \^etre vue
comme une mesure de la lissit\'e du morphisme pr\'ec\'edent.\\
La conjecture suivante a \'et\'e formul\'ee en 1989 \cite{Sp} (voir
aussi \cite{DS}) :\\

\textbf{Conjecture :} \textit{Toute fonction de Artin est born\'ee par
une fonction affine.}\\

Ce probl\`eme a \'et\'e \'etudi\'e, entre autres, par D. Delfino et I. Swanson, M. Hickel, S.
Izumi, M. Spivakovsky,  B. Teissier (cf. par exemple
\cite{Iz}, \cite{Sp}, ou \cite{DS}).\\
Voici quelques exemples d'id\'eaux dont la fonction de Artin est
born\'ee par une fonction affine (voir aussi \cite{Ro}) :\\
- les id\'eaux principaux de $\OO_N[X,\,Y,\,Z_1,...,\,Z_n]$ engendr\'es
par les polyn\^omes de la forme $XY-\sum_if_iZ_i$ o\`u l'id\'eal
$(f_1,...,\,f_n)$ est premier (th\'eor\`eme d'Izumi \cite{Iz}),\\
- les id\'eaux de $\OO_N[X_1,...,\,X_n]$ engendr\'es par des polyn\^omes
homog\`enes de degr\'e 1 (voir par exemple \cite{Ro}, th\'eor\`eme 3.1),\\
- les id\'eaux principaux de $\OO_N[X,\,X_1,...,\,X_n]$ engendr\'es par
les polyn\^omes de la forme $X^n+fX^{n-1}X_1+\cdots+f^{n}X_n$
o\`u $f\in\OO_N$ (\cite{DS}, th\'eor\`eme 3.10).\\
Enfin, pr\'ecisons que pour $N=1$, cette fonction apparait dans
l'\'etude de l'espace des arcs d'un germe singulier de vari\'et\'e
\cite{LJ} et en int\'egration motivique \cite{DL} (J. Denef et F. Loeser
relient cette fonction avec certaines s\'eries de Poincar\'e motiviques,
cf. th\'eor\`eme 8.1).\\

 Nous donnons ici un contre-exemple \`a cette conjecture :\\
\begin{theorem}\label{cont}
La fonction de Artin du polyn\^ome
$$P(X,\,Y,\,Z):=X^2-ZY^2\in\OO_N[X,\,Y,\,Z]$$
est born\'ee inf\'erieurement par une fonction polynomiale de degr\'e 2
si $N\geq 2$ et si la caract\'eristique de $\k$ est diff\'erente de 2.
\end{theorem}
Nous relions tout d'abord un r\'esultat de type approximation
diophantienne sur le corps $\K_N$ (le corps des fractions de $\OO_N$) \`a
l'existence d'une fonction de Artin (proposition \ref{homogene} et
corollaire \ref{approximation}). Nous montrons ensuite, \`a l'aide d'un
exemple (cf. exemple \ref{liouv}), qu'il n'existe pas de r\'esultat du
type th\'eor\`eme de Liouville pour une extension finie de $\K_N$ dans
$\wdh{\K}_N$, son compl\'et\'e pour la topologie $\m_N$-adique. Nous
nous appuyons alors sur cet exemple pour donner notre contre-exemple \`a
la conjecture cit\'ee pr\'ec\'edemment (th\'eor\`eme \ref{cont}). La
preuve du th\'eor\`eme \ref{cont} (partie \ref{preuve}) peut se lire
ind\'ependamment de la partie \ref{appdiop} qui traite de
l'approximation diophantienne. La partie \ref{appdiop} permet cependant
d'\'eclairer le th\'eor\`eme \ref{cont}.\\
 Enfin, nous d\'eduisons de cela, en derni\`ere partie, qu'il n'existe
pas de th\'eorie d'\'elimination des quantificateurs dans $\K_N$, pour
$N\geq 2$ (th\'eor\`eme \ref{elim}), pour le langage d\'efini dans la
derni\`ere partie.\\

Je suis tr\`es reconnaissant \`a M. Hickel et M. Spivakovsky pour leurs
conseils et encouragements. Je tiens aussi \`a remercier J. Nicaise pour
m'avoir fait remarquer que l'\'elimination des quantificateurs dans
$\K_N$ impliquait que toute fonction de Artin \'etait born\'ee par une
fonction affine.\\

Soient $N$ un entier positif non nul et $\k$ un corps ; nous utiliserons
les notations suivantes :
\begin{itemize}
\item
$\OO_N$ est l'anneau des s\'eries formelles $\Bbbk[[T_1,...,T_N]]$  et
$\m_N$ son id\'eal maximal (ou $\m$ quand il n'y aura pas d'\'equivoque
possible sur $N$).
\item $\ord$ est la valuation $\m_N$-adique sur $\OO_N$. Cette valuation
d\'efinit une norme $|\ |$ sur $\OO_N$ en posant $|x|=e^{-\ord(x)}$ et
cette norme induit une topologie appel\'ee topologie $\m$-adique.
\item $V_N:=\left\{\frac{x}{y}\,/\,x,\,y\in\OO_N \text{ et }
\ord(x)\geq\ord(y)\right\}$, l'anneau de valuation discr\`ete qui domine
$\OO_N$ pour $\ord$. Nous noterons $\m_N'$ son id\'eal maximal.
\item
 $\wdh{V}_N:=\k\left(\frac{T_1}{T_N},...,\,\frac{T_{N-1}}{T_N}\right)[[T_N]]$
est le compl\'et\'e  pour la topologie $\m_N'$-adique de $V_N$. En
effet, cet anneau correspond au compl\'et\'e le long du diviseur exceptionnel de l'anneau $\OO_N$. Nous noterons $\wdh{\m}_N$ l'id\'eal
maximal de cet anneau, $\ord$ l'extension de la valuation $\m_N'$-adique
et $|\ |$ l'extension de la norme associ\'ee.
\item $\K_N$ et $\wdh{\K}_N$ sont respectivement les corps de fractions
de $\OO_N$ et de $\wdh{V}_N$. On peut remarquer que $\wdh{\K}_N$ est le
compl\'et\'e de $\K_N$ pour la norme $|\ |$.
\end{itemize}
\begin{remark}
Un th\'eor\`eme de Wavrik (cf. \cite{Wa}) donne l'existence de la
fonction de Artin pour les syst\`emes d'\'equations polynomiales \`a
coefficients dans $\k\{T_1,...,\,T_N\}$, l'anneau des s\'eries
convergentes, quand $\k$ est muni d'une norme. Tout ce qui est
\'enonc\'e dans la suite est encore valable dans ce cas.
\end{remark}
\begin{remark}\label{Loj}
Soit $I=(f_1,...,\,f_p)$ un id\'eal de $\OO_N[X_1,...,\,X_n]$. Notons $F$
le sous ensemble de $\OO_N^n$ d\'efini par $f_1=\cdots=f_p=0$. Notons $d$
la distance sur $\OO_N^n$ induite par la norme $|\ |_{\infty}$ sur
$\OO_N^n$, elle-m\^eme d\'efinie par
$|(x_1,...,\,x_n)|_{\infty}=\max_i|x_i|$. Alors la fonction de Artin de
$I$ est born\'ee par une fonction affine si et seulement si nous avons
l'in\'egalit\'e suivante de type \L ojasiewicz :
$$\max_i\left|f_i(x_1,...,x_n)\right|\geq K\, d(x,\,F)^a$$
o\`u $K$ et $a$ sont des constantes ind\'ependantes des $x_i$ (pour une
preuve de ce r\'esultat, voir \cite{H2}).
\end{remark}

\section{Polyn\^ome homog\`ene \`a z\'ero isol\'e et approximation
diophantienne}\label{appdiop}

Nous allons montrer ici que la fonction de Artin d'un polyn\^ome
homog\`ene \`a z\'ero isol\'e est born\'ee par une fonction affine si et
seulement si nous avons un r\'esultat de type  approximation
diophantienne (cf. remarque \ref{rque}) :
\begin{proposition}\label{homogene}
Soit $P\in\OO_N[X_1,...,\,X_n]$ un polyn\^ome homog\`ene de degr\'e $d$
qui a pour unique z\'ero dans $\OO_N^n$ le point $(0,...,\,0)$. Notons
$Q_i(X_1,...,\,\wdh{X}_i,...,\,X_n)$ le polyn\^ome
$P(X_1,...,\,1,\,...,\,X_n)$ o\`u la variable $X_i$ est remplac\'ee par
1 dans $P$, pour $i\in\{1,...,\,n\}$.\\
Le polyn\^ome $P$ admet une fonction de Artin born\'ee par une fonction
affine  si et seulement si il existe deux constantes positives  $a$ et
$b$ telles que
nous ayons
\begin{equation}\label{liouville}\min_j\left\{\ord\left(\frac{u_j}{v}-y_j\right)\right\}\leq
a\,\ord(v)+b.\end{equation}
pour tout $i\in\{1,...,\,n\}$, pour toute racine
$(y_1,...,\,\wdh{y}_i,\,...,\,y_n)$ de $Q_i$ dans $\wdh{V}_N$ et pour
tout $u_1,...,\,\wdh{u}_i,\,...,\,u_n$ et $v$ dans $\OO_N$.
\end{proposition}

\begin{remark}\label{rque}
La condition (\ref{liouville}) peut aussi s'\'enoncer sous la forme
suivante : il existe deux constantes positives  $c$ et $K$ telles que
pour tout entier $i\in\{1,...,\,n\}$, pour toute racine
$(y_1,...,\,\wdh{y}_i,\,...,\,y_n)$ de $Q_i$ dans $\wdh{V}_N$ et pour
tout $u_1,...,\,\wdh{u}_i,\,...,\,u_n$ et $v$ dans $\OO_N$ nous ayons
$$\max_j\left\{\left|\frac{u_j}{v}-y_j\right|\right\}\geq K|v|^c.$$
\end{remark}
\begin{remark}\label{borneffect}
En fait, comme le montre la preuve de cette proposition, deux cas
peuvent se produire : soit aucun des polyn\^omes $Q_i$ n'admet de
solutions dans $\wdh{V}_N$, soit au moins un de ces polyn\^omes en admet
une. Dans le premier cas, la preuve nous permet de dire que la fonction
de Artin de $P$ est born\'ee par une fonction de la forme $i\lgm di+c$
o\`u $c\in\N$.\\
Dans le second cas, si la condition (\ref{liouville}) est v\'erifi\'ee,
alors  la fonction de Artin de $P$ est born\'ee par une fonction de la
forme $i\lgm d'i+c$ o\`u $c\in\R_{+}$ et $d'>d$ sont des constantes.
\end{remark}
Avant de donner la preuve de cette proposition, nous allons d'abord
\'enoncer un lemme qui va nous permettre de reformuler le th\'eor\`eme
de Artin dans le cas particulier des polyn\^omes \`a z\'ero isol\'e :\\
\begin{lemme}\label{ineg}
Soit $P$ un polyn\^ome de $\OO_N[X_1,...,\,X_n]$ ayant pour seul z\'ero
\` a coordonn\'ees dans $\OO_N$ le point $(0,...,\,0)$. Le polyn\^ome $P$
admet une fonction de Artin major\'ee par la fonction $\b$, si et
seulement si nous avons l'in\'egalit\'e :\\
\begin{equation}\label{ine1}\forall (x_1,...,\,x_n)\in\OO_N^n ,\
\ord(P(x_1,...,\,x_n))\leq \b\left(\min_k\{\ord(x_k)\}\right).\end{equation}
En particulier $P$ admet une fonction de Artin major\'ee par une
fonction affine, si et seulement si il existe deux constantes $a$ et $b$
telles que
\begin{equation}\label{ine}\forall (x_1,...,\,x_n)\in\OO_N^n ,\
\ord(P(x_1,...,\,x_n))\leq a\min_k\{\ord(x_k)\}+b.\end{equation}
\end{lemme}
\textit{D\'emonstration du lemme. }
Supposons que la fonction de Artin de $P$ soit  major\'ee par la
fonction $\b$. Soient $(x_1,...,\,x_n)\in\OO_N^n$ et $i\in\N$ tels que
$P(x_1,...,\,x_n)\in \m^{\b(i)+1}$ et $P(x_1,...,\,x_n)\notin
\m^{\b(i+1)+1}$. Alors d'apr\`es le th\'eor\`eme \ref{Ar}, il existe
$(\ovl{x}_1,...,\,\ovl{x}_n)\in\OO_N^n$ tel que pour tout $k$,
$x_k-\ovl{x}_k\in\m^{i+1}$, et $P(\ovl{x}_1,...,\,\ovl{x}_n)=0$. Par
hypoth\`ese sur $P$, nous avons n\'ecessairement $\ovl{x}_k=0$ pour tout
$k$. Et donc $\ord(x_k)\geq i+1$. Nous avons donc l'in\'egalit\'e
(\ref{ine1}).\\
Inversement si nous avons l'in\'egalit\'e (\ref{ine1}), et si
$(x_1,...,\,x_n)\in\OO_N^n$ v\'erifie $P(x_1,...,\,x_n)\in\m^{\b(i)+1}$,
alors n\'ecessairement $x_k\in\m^{i+1}$ pour tout $k$, et la fonction de
Artin de $P$ est born\'ee par la fonction $\b$.$\quad\Box$\\
\\
\textit{D\'emonstration de la proposition. }
Soient $P$ comme dans l'\'enonc\'e et $x_1,...,\,x_n\in\OO_N$. Quitte \`a
changer le nom des variables, nous pouvons supposer que nous avons
$\ord(x_1)\leq...\leq\ord(x_n)$. Nous avons alors
$$P(x_1,...,\,x_n)=x_1^dP\left(1,\frac{x_2}{x_1},...,\,\frac{x_n}{x_1}\right)
.$$
Notons alors $Q_1(Y_2,...,\,Y_n)$ le polyn\^ome $P(1,\,Y_2,...,\,Y_n)$.
Le polyn\^ome $Q_1$ n'a pas de racine dans $\OO_N$. Il n'en a donc pas
non plus dans $V_N$, mais peut en avoir dans $\wdh{V}_N$.\\
\\
Supposons que $Q_1$ n'ait pas de racine dans $\wdh{V}_N$. Comme
$\wdh{V}_N$ est de la forme $\K[[T]]$ o\`u $\K$ est un corps et $T$ une
variable formelle, d'apr\`es le th\'eor\`eme de Greenberg (cf.
\cite{Gr}), $Q_1$ admet une fonction de Artin born\'ee par une constante
$c$ et dans ce cas nous obtenons 
 $$\ord\left(P\left(1,\frac{x_2}{x_1},...,\,\frac{x_n}{x_1}\right)\right)<
c+1.$$ Donc
$$\ord\left(P(x_1,...,\,x_n)\right)<d\,\ord(x_1)+c+1\ .$$
Et donc, d'apr\`es le lemme \ref{ineg}, la fonction de Artin de $P$ est
born\'ee par une fonction affine dont le coefficient de lin\'earit\'e
peut \^etre choisi \'egal \`a $d$. Ceci correspond au premier cas
\'enonc\'e dans la remarque \ref{borneffect}.\\
\\
Si $Q_1$ a des racines dans $\wdh{V}_N$, toujours d'apr\`es \cite{Gr},
$Q_1$ admet  une fonction de Artin born\'ee par une fonction affine
$i\lgm \la i+\mu$. Notons $(y_2,...,\,y_n)$ un plus proche z\'ero de
$\left(\frac{x_2}{x_1},...,\,\frac{x_n}{x_1}\right)$ pour la topologie
$\m$-adique.\\
Supposons que nous ayons
$\min_j\left\{\ord\left(\frac{x_j}{x_1}-y_j\right)\right\}\leq
a\,\ord(x_1)+b$ o\`u $a$ et $b$ sont des constantes. Alors nous avons
$$\ord\left(Q_1\left(\frac{x_2}{x_1},...,\,\frac{x_n}{x_1}\right)\right)<\la\left(\min_j\left\{\ord\left(\frac{x_j}{x_1}-y_j\right)\right\}+1\right)+\mu+1$$
$$\qquad\qquad\quad\ <\la(a\,\ord(x_1)+b+1)+\mu+1.$$
D'o\`u
$$\ord(P(x_1,...,\,x_n))<(a\la+d)\,\ord(x_1)+\la a(b+1)+\mu+1.$$
Et donc, l\`a encore d'apr\`es le lemme \ref{ineg}, la fonction de Artin
de $P$ est born\'ee par une fonction affine dont le coefficient de
lin\'earit\'e est strictement plus grand que $d$. Ceci correspond au
second cas \'enonc\'e dans la remarque \ref{borneffect}.\\
\\
Inversement, supposons qu'il existe $i\in\{1,...,\,n\}$ tel que, pour
tout $c\in\N$, il existe une racine
$(y_{1,c},...,\,\wdh{y}_{i,c},\,...,\,y_{n,c})$ de $Q_i$ dans
$\wdh{V}_N$ et des $u_{j,c}$ et $v_c$ dans $\OO_N$ tels que
$$\min_j\left\{\ord\left(\frac{u_{j,c}}{v_c}-y_{j,c}\right)\right\}\geq
c\,\ord(v_c).$$
Comme $\ord(y_{j,c})\geq 0$, nous avons
$\ord\left(\frac{u_{j,c}}{v_c}\right)\geq 0$ pour tout $j$ et pour tout
$c$ et donc
$\min(\ord(v_c),\,\ord(u_{2,c}),...,\,\ord(u_{n,c}))=\ord(v_c)$.\\
Nous avons
 $$Q_i\left(\frac{u_{1,c}}{v_c},...,\,\frac{u_{n,c}}{v_c}\right)=Q_i\left(\frac{u_{1,c}}{v_c},...,\,\frac{u_{n,c}}{v_c}\right)-Q_i\left(y_{1,c},...,\,y_{n,c}\right)=$$
 $$=Q_i\left(\frac{u_{1,c}}{v_c},...,\,\frac{u_{n,c}}{v_c}\right)-Q_i\left(y_{1,c},\,\frac{u_{2,c}}{v_c},...,\,\frac{u_{n,c}}{v_c}\right)+$$
$$+Q_i\left(y_{1,c},\,\frac{u_{2,c}}{v_c},...,\,\frac{u_{n,c}}{v_c}\right)-\cdots-Q_i\left(y_{1,c},...,\,y_{n,c}\right)$$
Or nous avons
$$Q_i\left(y_{1,c},...,y_{j-1,c},\frac{u_{j,c}}{v_c},...,\frac{u_{n,c}}{v_c}\right)-Q_i\left(y_{1,c},...,y_{j,c},\frac{u_{j+1,c}}{v_c},...,\frac{u_{n,c}}{v_c}\right)=$$
$$=\left(\frac{u_{j,c}}{v_c}-y_{j,c}\right)U_i$$
o\`u $U_i\in \wdh{V}_N$ (car $Q_i$ est \`a coefficients dans $\OO_N$, et
les $\frac{u_{j,c}}{v_c}$ et les $y_{j,c}$ sont dans $\wdh{V}_N$).\\
Nous en d\'eduisons que
\begin{equation}\label{inegalitecle}\ord\left(Q_i\left(\frac{u_{2,c}}{v_c},...,\,\frac{u_{n,c}}{v_c}\right)\right)\geq
\min_j\left\{\ord\left(\frac{u_{j,c}}{v_c}-y_{j,c}\right)\right\}\geq
c\,\ord(v_c).\end{equation}
D'o\`u
$$\ord(P(v_c,\,u_{2,c},...,\,u_{n,c}))\geq (c+d)\,\ord(v_c).$$
Donc $P$ n'admet pas de fonction de Artin born\'ee par une fonction
affine, car $P$ ne v\'erifie pas l'in\'egalit\'e (\ref{ine}) du lemme
\ref{ineg}.$\quad\Box$\\
\\
Nous d\'eduisons de la proposition pr\'ec\'edente le r\'esultat
d'approximation diophantienne suivant :
\begin{proposition}\label{approximation}
Soit $x\in\wdh{\K}_N$ alg\'ebrique sur $\K_N$ tel que $x\notin\K_N$.
Alors il existe une fonction $\g :\N\lgw \N$ telle que
$$\ord\left(\frac{u}{v}-x\right) \leq \g(\ord(v))$$
pour tous $u$ et $v$ dans $\OO_N$.\\
Si $Q=a_dX^d+\cdots+a_0$ est le polyn\^ome minimal de $x$ sur $\OO_N$,
alors $\g$ est major\'ee par une fonction affine si et seulement  si la
fonction de Artin de
$P=a_dX^d+a_{d-1}X^{d-1}Y\cdots+a_0Y^d\in\OO_N[X,\,Y]$ est major\'ee par
une fonction affine, c'est-\`a-dire qu'il existe $a\geq 1$ et $K\geq 0$
tels que
$$\left|\frac{u}{v}-x\right|\geq K|v|^{a},\ \forall u,\,v\in\OO_N.$$
\end{proposition}
\begin{proof} Soit  $x\in\wdh{\K}_N$ alg\'ebrique sur $\K_N$. Si
$\ord(x)\geq 0$, soit $Q$ un polyn\^ome irr\'eductible de $\OO_N[X]$ tel
que $Q(x)=0$ et $P(X,Y)$ l'homog\'en\'eis\'e de $Q$. Le corollaire
d\'ecoule alors de la proposition pr\'ec\'edente si $P$ admet une
fonction de Artin born\'ee par une fonction affine. Si ce n'est pas le
cas, l'existence de la fonction de Artin de $P$, not\'ee $\b$, nous
permet de dire, en reprenant l'in\'egalit\'e (\ref{inegalitecle}) de la
derni\`ere partie de la preuve de la proposition pr\'ec\'edente et grace
au lemme \ref{ineg},
$$\b(\ord(v))\geq\b(\min\{\ord(u),\ord(v)\})\geq\ord(P(u,v))=$$
$$=\ord\left(Q\left(\frac{u}{v}\right)\right)+d\,\ord(v)\geq\ord\left(\frac{u}{v}-x\right)+d\,\ord(v).$$
Le r\'esultat en d\'ecoule directement.\\
Si $\ord(x)<0$, notons $Q$ un polyn\^ome irr\'eductible de $\OO_N[X]$ tel
que $Q(x)=0$. Soit $P(X,\,Y)$ l'homog\'en\'eis\'e de $Q$ et
$R(Y)=P(1,\,Y)$. Le polyn\^ome $R$ est irr\'eductible, de m\^eme degr\'e
que $Q$ et $R(1/x)=0$. Nous avons alors
$Q\left(\frac{u}{v}\right)=\frac{u^d}{v^d}R\left(\frac{v}{u}\right)$. Si
le terme initial de $\frac{u}{v}$ est diff\'erent du terme initial de
$x$ pour la valuation $\ord$, alors nous obtenons
$\ord\left(\frac{u}{v}-x\right) \leq \ord(x)$. Sinon,
$\ord(u)+\ord(x)=\ord(v)$. Sachant que $1/x\in\wdh{V}_N$, d'apr\`es le
lemme \ref{ineg} et l'existence la fonction de Artin $\b$ du polyn\^ome
$P$, nous obtenons comme pr\'ec\'edemment l'in\'egalit\'e
$\ord\left(\frac{v}{u}-\frac{1}{x}\right) \leq
\b(\ord(u))-d\,\ord(u)\leq \b(\ord(u))$.\\ Donc nous avons
$$\ord\left(\frac{u}{v}-x\right)=\ord\left(\frac{ux}{v}\left(\frac{v}{u}-\frac{1}{x}\right)\right)=\ord\left(\frac{v}{u}-\frac{1}{x}\right)
\leq \b(\ord(v)-\ord(x)).$$ Dans tous les cas, nous avons
$\ord\left(\frac{u}{v}-x\right)\leq \b(\ord(v)-\ord(x))+\ord(x)$.\\
\end{proof}
\begin{example}\label{liouv}
Supposons que la caract\'eristique de $\k$ est diff\'erente de 2. Nous
pr\'esentons ici une suite $(x_p)_{p\in\N\backslash{\{0,\,1,\,2\}}}$
d'\'el\'ements de $\wdh{\K}_N$, chacun de degr\'e 2 sur $\K_N$, pour
lesquels il existe $u_{p,\,k}$ et $v_{p,\,k}$ tels que
$\left|x_p-\frac{u_{k,\,p}}{v_{k,\,p}}\right|=C_p|v_k|^{\frac{p}{2}-1}$
et $\ord(v_{p,\,k})$ tend vers $+\infty$ avec $k$, o\`u $C_p$ est une
constante qui d\'epend de $p$. Nous voyons donc que s'il existe une
version lin\'eaire d'approximation diophantienne pour les s\'eries en
plusieurs variables, contrairement au cas des nombres r\'eels
alg\'ebriques, la meilleure borne $c$ telle qu'il existe $K$ avec
$$\left|x-\frac{u}{v}\right|>K|v|^c,\ \forall u,\,v\in\OO_N$$
ne peut pas \^etre born\'ee par le degr\'e de l'extension de $\K_N$ par
$x$. Il n'existe donc pas de version du th\'eor\`eme de Liouville pour
les extensions finies de $\K_N$ dans $\wdh{\K}_N$.\\
\\
Soit $P_p(X,Y)$ le polyn\^ome $X^2-(T_1^2+T_2^p)Y^2$ avec
$p\in\N\backslash\{0,\,1\,\,2\}$. Le terme $T_1^2+T_2^p$ n'est pas un
carr\'e car la caract\'eristique du corps de base est diff\'erente de 2.
Le seul z\'ero de $P_p$ est alors $(0,\,0)$. Soit
$Q_p(X)=X^2-(T_1^2+T_2^p)$. Le polyn\^ome $Q_p$ a deux racines dans
$\wdh{V}_N$ qui sont $x_p$ et $-x_p$ avec
$$x_p=
T_1\left(1+\frac{1}{2}\frac{T_2^p}{T_1^2}-\frac{1}{8}\frac{T_2^{2p}}{T_1^4}+\cdots+\frac{(-1)^{n-1}(2n-2)!}{2^{2n-1}(n-1)!n!}\frac{T_2^{np}}{T_1^{2n}}+\cdots\right)$$
Notons $a_n:=(-1)^{n-1}\frac{(2n-2)!}{2^{2n-1}(n-1)!n!}$.
Soit $k$ un entier positif et notons
$$x_{p,\,k}:=T_1\sum_{i=0}^{k-1}a_i\frac{T_2^{ip}}{T_1^{2i}}.$$
Nous avons
\begin{equation}\label{cle}x_p-x_{p,\,k}=T_1\sum_{i\geq
k}a_i\frac{T_2^{ip}}{T_1^{2i}}\in\wdh{\m}_N^{(p-2)k+1}\end{equation}
et $x_{p,\,k}$  s'\'ecrit sous la forme $\frac{u_{p,\,k}}{v_{p,\,k}}$
avec $u_{p,\,k}$ et $v_{p,\,k}$ premiers entre eux et
$v_{p,\,k}=T_1^{2k-3} $. D'o\`u
$$\ord\left(x_p-\frac{u_{p,\,k}}{v_{p,\,k}}\right)=k(p-2)+1=\left(\frac{p}{2}-1\right)\ord(v_{p,\,k})+\frac{3}{2}p-2,$$
$$\text{ ou encore }
\left|x_p-\frac{u_{p,\,k}}{v_{p,\,k}}\right|=e^{-\frac{3}{2}p+2}|v_{p,\,k}|^{\frac{p}{2}-1}$$
\end{example}
\begin{example}
Plus g\'en\'eralement, S. Izumi a \'etudi\'e les polyn\^omes de la forme
$X^d-aY^d$ o\`u $a$ n'est pas une puissance $d$-i\`eme dans $\OO_N$, et a
montr\'e que la fonction de Artin de ces polyn\^omes est born\'ee par
une fonction affine (cf. proposition 5.1 \cite{Iz}), sans pour autant
donner de bornes explicites.
\end{example}

\section{Preuve du th\'eor\`eme \ref{cont}}\label{preuve}
Nous pouvons  donner la preuve du th\'eor\`eme \ref{cont} annonc\'e dans
l'introduction. Celle-ci consiste \` a interpr\'eter diff\'eremment le
contre-exemple \ref{liouv} \`a l'approximation diophantienne de type
Liouville :\\
\\
Soit $N\geq 2$ fix\'e. Nous notons $P(X,\,Y,\,Z):=X^2-Y^2Z$. Supposons
la caract\'eristique de $\k$ diff\'erente de 2. Soient $p$ et $k$ des
entiers strictement plus grands que 2 et soient
$$u_{p,\,k}=T_1^{2k-2}\sum_{i=0}^{k-1}a_i\frac{T_2^{ip}}{T_1^{2i}},\
v_k=T_1^{2k-3},\text{ et }z_p=T_1^2+T_2^{p}$$
avec, pour tout $n\in\N$,
$a_n:=(-1)^{n-1}\frac{(2n-2)!}{2^{2n-1}(n-1)!n!}$.\\
Notons, comme dans l'exemple \ref{liouv},
$x_{p,\,k}:=\frac{u_{p,\,k}}{v_k}=T_1\sum_{i=0}^{k-1}a_i\frac{T_2^{ip}}{T_1^{2i}}$,
et $x_p:=T_1\sum_{i\geq 0}a_i\frac{T_2^{ip}}{T_1^{2i}}$. En particulier
$x_p^2=z_p$.\\
Alors
$$P(u_{p,\,k},\,v_k,\,z_p)=\left(T_1^2\left(\sum_{i=0}^{k-1}a_i\frac{T_2^{ip}}{T_1^{2i}}\right)^2-(T_1^2+T_2^{p})\right)T_1^{4k-6}=$$
$$=\left(\left(\frac{u_{p,\,k}}{v_k}\right)^2-z_p\right)v_k^2=\left(x_{p,k}-x_p\right)\left(x_{p,\,k}+x_p\right)v_k^2\in\wdh{\m}_N^{(p+2)k-4}$$
d'apr\`es la relation (\ref{cle}) de l'exemple \ref{liouv}.\\
\\
Si $(x,\,y,\,z)$ est un z\'ero de $P$ alors soit $z$ est un carr\'e,
soit $x=y=0$. Or
$$\sup_{t\in\OO_N}(\ord(z_p-t^2))= p$$
car la caract\'eristique de $\k$ est diff\'erente de 2, et
$$ \min(\ord(u_{p,\,k}),\,\ord(v_k))=2k-3.$$
Donc, en posant $p=k-2$, nous avons
$$P(u_{k-2,\,k},\,v_{k},\,z_{k-2})\in\m^{k^2-4}\ $$
$$\text{et }\sup\big(\min\{\ord(u_{k-2,\,k}-x),
\,\ord(v_{k}-y),\,\ord(z_{k-2}-z)\}\big)\leq 2k-3$$
o\`u la borne sup\'erieure est prise sur toutes les racines
$(x,\,y,\,z)$ de $P$.\\
\\
Notons $\b_N$ la fonction de Artin de $P$. Soient $i\in\N^*$ pair et
$k\in\N$ tel que $i=2k-2$. D'apr\`es ce qui pr\'ec\`ede, il existe une
solution approch\'ee de $P$ \`a l'ordre
$\left(\frac{i+2}{2}\right)^2-4=\left(\frac{i+2}{2}\right)^2-5+1$, mais
il n'existe aucune solution de $P$ "proche" de cette solution \`a
l'ordre $i+1$. Donc $\b_N(i)\geq \left(\frac{i+2}{2}\right)^2-5$. Si $i$
est impair, comme $\b_N(i)\geq\b_N(i-1)$, nous avons
$\b_N(i)\geq\left(\frac{i}{2}\right)^2-5$. Donc pour tout $i\in\N^*$,
nous avons
$$\b_N(i)\geq \frac{i^2}{4}-5\, .\qquad\qquad\Box$$\\
\section{Application : non-existence d'\'elimination des quantificateurs
dans le corps $\k((T_1,...,\,T_N))$ pour $N\geq 2$}
J. Denef et F. Loeser \cite{DL} ont donn\'e une preuve du th\'eor\`eme
de Greenberg \`a l'aide d'un r\'esultat d'\'elimination des
quantificateurs d\^u \`a J. Pas \cite{Pa}. Nous pouvons appliquer ici la
m\^eme m\'ethode pour montrer que dans le cas $N\geq 2$, il n'y a pas
d'existence d'une th\'eorie d'\'elimination des quantificateurs dans le
corps $\k((T_1,...,\,T_N))$ muni du langage d\'efini ci-dessous.\\
Dans  cette partie,  $\k$ est un corps alg\'ebriquement clos.\\
Soit $\LL_{Pre}$ le langage du premier ordre dont les variables sont les
\'el\'ements de $\Z$ et les symboles sont $+$, $\leq$, $0$, $1$ et pour
tout $d\in\Z_{\geq2}$  un symbole pour signifier la relation binaire
$x\equiv y$ mod $d$.\\
Soit $\LL_{\k}$ le langage du premier ordre dont les variables sont les
\'el\'ements de $\k$ et les symboles sont  $+$, $-$, $\times$, $0$, $1$.\\
Soit $\LL_N$ le langage du premier ordre dont les variables sont les
\'el\'ements de $\K_N$ et les symboles sont  $+$, $-$, $\times$, $0$,
$1$ et le symbole pour la fonction $\ord$.\\
Nous pouvons alors \'enoncer le r\'esultat suivant d\^u \`a F. Delon
\cite{D} :
\begin{theorem}\label{elim}
Soit $N\geq2$ fix\'e. Consid\'erons le langage du premier ordre
\`a trois sortes $\LL:=(\LL_{Pre},\,\LL_{\k},\,\LL_N,\,\ord,\,\pi)$, o\`u
$\pi$ est une fonction de $\K_N$ vers $\k$. Soit $\LL'$ le langage
form\'e \`a partir de $\LL$ en lui ajoutant autant de symboles que
l'on veut de telle sorte que la restriction de $\LL'$ \`a $\Z$ soit
\'egale \`a $\LL_{Pre}$. Alors $\LL'$ n'admet pas d'\'elimination des
$\K_N$-quantificateurs.\\
\end{theorem}
\begin{proof} Consid\'erons $\b$, la fonction de Artin du polyn\^ome
$P=X^2-ZY^2$ vu comme polyn\^ome de $\OO_N[X,\,Y,\,Z]$. Cette fonction
peut se d\'efinir dans les langages $\LL$ et $\LL'$, et donc son graphe
aussi. En effet, pour tout entier $n$ nous avons
$$\left( P(x,\,y,\,z) \in\m^{\b(n)+1}\Longrightarrow\left( \exists
\ovl{x},\,\ovl{y},\,\ovl{z}\,/\,\left(
x-\ovl{x},\,y-\ovl{y},\,z-\ovl{z}\in\m^{n+1}\right)\right.\right.$$
$$\left.\left.\wedge P(\ovl{x},\,\ovl{y},\,\ovl{z})=0\right)\right)$$
$$\wedge\left(P(x,\,y,\,z) \in\m^{\b(n)}\Longrightarrow\left(
\neg\exists \ovl{x},\,\ovl{y},\,\ovl{z}\,/\,\left(
x-\ovl{x},\,y-\ovl{y},\,z-\ovl{z}\in\m^{n+1}\right)\right.\right.$$
$$\left.\left.\wedge P(\ovl{x},\,\ovl{y},\,\ovl{z})=0\right)\right)$$
Si $\LL'$ admettait une \'elimination des $\K_N$-quantificateurs, alors
d'apr\`es un r\'esultat de Presburger, qui dit que $\LL_{Pre}$ admet une
\'elimination des quantificateurs (cf. \cite{Pr}), et d'apr\`es le
th\'eor\`eme de constructibilit\'e de Chevalley, le langage $\LL'$
admettrait une \'elimination des quantificateurs. Le graphe de $\b$,
inclus dans $\Z^2$, serait donc semi-alg\'ebrique dans le langage
$\LL_{Pre}$ et il existerait alors une partition finie de $\N$ en
classes de congruences telle que $\b$ soit affine sur chacune de ces
classes, ce qui est faux d'apr\`es ce qui pr\'ec\`ede.
\end{proof}
\begin{remark}
Pour $N=1$ un r\'esultat d'\'elimination des quantificateurs a \'et\'e
obtenu par J. Pas (cf. \cite{Pa}).
\end{remark}

\end{document}